%zudilin.tex: 
%%a Plain TeX file by Robert Dougherty-Bliss and Doron Zeilberger (x pages)

%begin macros

\baselineskip=14pt
\parskip=10pt

\magnification=\magstephalf

\def\1{{\overline{1}}}
\def\2{{\overline{2}}}
\parindent=0pt
\overfullrule=0in

\def\frac#1#2{{#1 \over #2}}
%\headline={\rm  \ifodd\pageno  \RightHead  \else  \LeftHead  \fi}
%\def\RightHead{\centerline{
%Title
%}}
%\def\LeftHead{ \centerline{Doron Zeilberger}}
%end macros
\centerline
{\bf Exploring General Ap\'ery Limits via the Zudilin-Straub t-transform
 }
\bigskip
\centerline
{\it Robert DOUGHERTY-BLISS and Doron ZEILBERGER}

{\bf Abstract:} Inspired by a recent beautiful construction of Armin Straub and Wadim Zudilin, that `tweaked' the sum of the $s^{th}$ powers of
the $n$-th row of Pascal's triangle, getting instead of sequences of numbers,  sequences of rational functions,
we do the same for general binomial coefficients sums, getting a practically unlimited supply of Ap\'ery limits. While getting
what we call ``major Ap\'ery miracles", proving irrationality of the associated constants (i.e. the so-called Ap\'ery limits)
is very rare, we do get, every time, at least a   ``minor Ap\'ery miracle" where  an explicit constant, defined as an (extremely slowly-converging)
limit of some explicit sequence,  is expressed as an Ap\'ery limit of some recurrence, with some initial conditions, thus enabling a very fast computation of that
constant, with exponentially decaying error.

{\bf Preface: The Major and Minor (but still interesting!)  Ap\`ery  Miracles}

One way that Roger Ap\'ery's [A] [vdP] seminal proof of the irrationality  of $\zeta(3)$  could have been discovered, in a {\it counterfactual world}, was to
consider, {\it out of the blue}, the second-order linear recurrence
$$
n^3 u_n \,- \, \left(17 n^{2}+51 n +39\right) \left(2 n +3\right)\,u_{n-1} \,+ \, (n-1)^3 u_{n-2}=0 \quad,
$$
and let $a_n$ and $b_n$ be the solutions of that recurrence with {\bf initial conditions}
$$
a_0=0, a_1=6 \quad; \quad
b_0=1, b_1=5 ,
$$
then let the computer compute many terms, evaluate $\frac{a_{1000}}{b_{1000}}$ to many decimals, and then use Maple's {\tt identify}, and {\it lo and behold}, get that it (most probably) equals $\zeta(3)$
(i.e. $\sum_{i=1}^{\infty} \frac{1}{i^3}$). Then, still {\it empirically} and {\it numerically}, after rewriting $\frac{a_n}{b_n}$ as $\frac{a'_n}{b'_n}$, where now
{\bf both} numerator and denominators are integers (initially $b_n$ were integers, but $a_n$ were not), estimate that there exists a {\it positive} number $\delta$ (about $0.0805$) such that
$$
|\frac{a'_n}{b'_n} - \zeta(3)| \leq \frac{CONSTANT}{(b'_n)^{1+\delta}} \quad ,
$$
that immediately entails (see [vdP]) that $\zeta(3)$ is irrational.

Using the terminology that has now become standard (e.g. the titles of [CS] and [SZ]), we say that $\zeta(3)$ is the {\it Ap\'ery limit} of the above recurrence, and initial conditions.

The reason that this was a {\bf major} miracle, as explained so eloquently in [vdP], is that while {\bf any} naturally occurring constant that is not {\bf obviously} rational,
e.g. the sum of the series $\sum_{i=1}^{\infty} \frac{1}{i(i+1)}$, is {\bf definitely} (in the everyday sense of the word) irrational (there are only $\aleph_0$ rational numbers,
while there are $2^{\aleph_0}$ real numbers, hence the `probability' of a real number being rational is a (very small!) $0$), it is extremely difficult to
(rigorously) {\bf prove} that a {\it specific} constant is irrational. Witness the fact that, in spite of many attempts,  there are still no proofs of
the irrationality of the Euler-Mascheroni constant $\gamma$ (the limit of the partial sums of the harmonic series minus $\log(n)$, or equivalently $-\int_0^{\infty} e^{-x} \log x \, dx$),
the Catalan constant $C:=\sum_{i=0}^{\infty} \frac{(-1)^i}{(2i+1)^2}$, $\zeta(5):=\sum_{i=1}^{\infty} \frac{1}{i^5}$ (and more generally $\zeta(2i+1)$, for all $i \geq 2$).
The same is also true for $e+\pi$ and $e\,\pi$ (but they can't be {\bf both} rational [proof left to the reader]).

This theme is pursued in [Ze2], [Ze3] and much more recently, in [DKZ], [DZ] and [ZeZu], where the motivation was to discover {\bf irrationality proofs} of other constants.
In [DKZ] there were quite a few `Ap\'ery miracles', alas, most of them were reproved irrationality of constants that were already proved irrational (for example, algebraic numbers or logarithms of them) and the novelty was establishing
{\it explicit  irrationality measures}. This is  still interesting, but not exciting.
We also  found  few other `weird' constants given in terms of products of Gamma values at
rational numbers, that should yield fully rigorous {\bf first} proofs of {\bf explicit} constants, but since these constants were name-less, it gave us neither fame nor fortune.

But the {\bf minor} Ap\'ery miracle was not {\bf number-theoretical} but rather {\bf numerical-analytical}. Here is an {\bf explicitly} defined constant
$\zeta(3):=\sum_{i=1}^{\infty} \frac{1}{i^3}$, in other words the limit of the sequence of rational numbers
$\{\sum_{i=1}^{n} \frac{1}{i^3}\}$, that converges {\bf very slowly} to its limit.
Realizing it as an {\it Ap\'ery limit}, i.e.
coming up with an explicit {\bf linear recurrence equation with polynomial coefficients},
and two sets of {\bf initial conditions}, for which the limit of the ratios (of the emerging two sequences) converge to that number with an {\bf exponentially decaying error}.
This enables one to compute the constant in question to many decimal digits. On the other hand, computing it to that accuracy using the
definition would take zillions of years. The purpose of this article is to show how one can produce lots of other `minor Ap\'ery miracles'
where one can express {\it explicitly} defined constants, whose definition entails very slow convergence, as Ap\`ery limits of recurrences and initial conditions
that enables computing these constants with exponentially decaying errors. The key idea is to introduce what we will call the {\bf Zudilin-Straub t-Transform},
that generalizes a construction Armin Straub and Wadim Zudilin discovered [SZ] for the special case of {\bf sums of powers of binomial coefficients}.

But before defining the Zudilin-Straub t-transform, let's define  {\it Ap\'ery limit} more formally, and also introduce the new notion of  {\it Generalized Ap\'ery limit}.

{\bf Definition 1}: A constant $c$ is an {\it Ap\'ery limit} if there exists a  {\bf homogeneous} linear recurrence with {\bf polynomial coefficients}
$$
\sum_{i=0}^{L} p_i(n)\, X(n+i) \, = \, 0 \quad,
$$
and two sets of {\bf initial conditions} $[a_0, a_1, \dots, a_{L-1}]$ and  $[b_0, b_1, \dots, b_{L-1}]$, such that if $A(n)$ and $B(n)$ are the solutions of that same recurrence
with 
$$
A(0)=a_0 \quad, \dots, \quad A(L-1)=a_{L-1} \quad ; \quad
B(0)=b_0 \quad, \dots, \quad B(L-1)=b_{L-1} \quad ,
$$
then
$$
c= \, \lim_{n \rightarrow \infty} \frac{A(n)}{B(n)} \quad .
$$

While not part of the definition, it turns out that often, and in all the naturally occurring cases, we also have the following nice feature.

{\bf Exponential decay of Error Property}

There exist  constants $C>0$ and $\alpha>1$ such that
$$
|\frac{A(n)}{B(n)} \, - \, c | \, \leq \, \frac{C}{\alpha^n} \quad .
$$

Given a recurrence and initial conditions, it is very fast to compute many terms. In fact one only needs {\bf constant memory} (well, linear memory if you go by bit-size) and {\bf linear time} to compute any specific approximation, $\frac{A(n)}{B(n)}$.

Let's introduce a mild extension, that of a  {\it Generalized Ap\'ery limit}.

{\bf Definition 2}: A constant $c$ is a {\it Generalized Ap\'ery limit} if there exist two sequences of rational numbers $A(n)$ and $B(n)$ such that
$$
c= \, \lim_{n \rightarrow \infty} \frac{A(n)}{B(n)} \quad ,
$$
where $A(n)$ and $B(n)$ are solutions of linear recurrences with polynomial coefficients (the first homogeneous, the second inhomogeneous)
$$
\sum_{i=0}^{L} p_i(n)\, B(n+i) \, = \, 0 \quad,
$$
subject to initial conditions $B(0)=b_0, \dots, B(L-1)=b_{L-1}$, and
$$
\sum_{i=0}^{L} p_i(n)\, A(n+i) \, = \, C(n) \quad,
$$
subject to initial conditions $A(0)=a_0, \dots, A(L-1)=a_{L-1}$, where the {\bf right side}, $C(n)$, in turn is a solution of another, (this time homogeneous) linear recurrence equation
with polynomial coefficients
$$
\sum_{i=0}^{M} q_i(n)\, C(n+i) \, = \, 0 \quad,
$$
subject to some initial conditions $C(0)=c_0, \dots, C(M-1)=c_{M-1}$.

Note that by using recurrence operators, it is easy to express both $A(n)$ and $B(n)$ as solutions of the {\bf same} homogeneous linear recurrence, so a generalized Ap\'ery limit can always be expressed
as an Ap\'ery limit, alas, with the recurrences being of a much higher order, namely $L+M$. 

Also note that in order to prove irrationality {\it \`a la Ap\'ery}, exponential decay of error does not suffice. After writing the  quotients of rational numbers $A(n)/B(n)$ as
$\frac{a_n}{b_n}$, where $a_n$ and $b_n$ are {\bf integers}, one needs an inequality of the form,
$$
|\frac{a_n}{b_n} \, - \, c | \, \leq \, \frac{CONSTANT}{b_n^{1+\delta}} \quad ,
$$
where $\delta>0$, yielding an {\it irrational measure} $1+\frac{1}{\delta}$ (see [vdP]).
The attempts in [Ze2],[Ze3], [DZ], and [DKZ], to find new (major) Ap\'ery miracles, i.e. irrationality proofs of hopefully new constants,
consisted of going {\it backwards}. Rather than trying to hit the ``bull's-eye", one shoots first, and then draws the bull's eye around the bullet hole. 
Using the Zeilberger, Almkvist-Zeilberger, and multi-Almkvist Zeilberger algorithms (see [DKZ] for references) we generated
recurrences for known sequences $B(n)$, then changed the initial conditions, getting a companion sequence $A(n)$, then we computed (very fast, to high accuracy)
approximations to the limit of $A(n)/B(n)$. Then we used Maple's {\tt identify} (and our extensions) to {\bf conjecture} an explicit value of the
Ap\'ery limit, and hoped to prove it later. But often, neither Maple, nor our extension, was able to {\it identify} the Ap\'ery limit.
In many cases we were able to identify the Ap\'ery limit, but the $\delta$ turned out to be negative, so it was useless for proving irrationality.
Nevertheless, since we always had $1+\delta>0$, we still got {\bf exponentially decaying error}.

In this paper we will forget about our irrationality obsession and only enjoy the exponential decay of error property of Ap\'ery limits (and the generalized version).
However here the focus would be to introduce  {\bf explicit} constants, defined as a limit of very slowly-converging sequences, and express them as Ap\'ery limits
or generalized Ap\'ery limits, enabling computing these constants very fast, to any desired accuracy. In other words, we will do what numerical analysts call {\it convergence acceleration},
with {\bf very dramatic} acceleration.

{\bf The work of Straub and Zudilin that motivated the  Zudilin-Straub t-transform}

In a recent beautiful article [SZ] (Theorem 1.3 there) (that brilliantly proved some conjectures in the equally beautiful article [CS]) they expressed
$\zeta(2j)$ $j=1,2,3\dots$ (or equivalently $\pi^{2j}$) as Ap\'ery limits with {\bf explicit} recurrences and {\bf explicit} initial conditions. The recurrences in question were obtained from
the {\bf Zeilberger} [Ze1] (see also [PWZ]) algorithm applied to {\bf sums of binomial powers}, also known as {\bf Franel numbers} (here $s$ is any positive integer):
$$
F^{(s)}_n:=\sum_{k=0}^{n} \, {{n} \choose {k}}^s \, = \,
\sum_{k=0}^{n} \, \frac{n!^s}{k!^s (n-k)!^s} 
\quad,
$$
with $s \geq 2j+1$. The beauty is that they did it for {\it infinitely many cases}.
 
Let's describe what they did. Recall that the {\bf rising factorial}, also called {\bf Pochhammer symbol} (that features in the definition of a {\bf hypergeometric series}), is defined by
$$
(x)_n \, := \, x\,(x+1)\, \cdots \, (x+n-1) \quad .
$$
(Note that $(1)_n=n!$.) 

The starting point in [SZ] was to consider, instead of the sequence of {\bf integers} $\{F^{(s)}_n\}_{n=1}^{\infty}$, the sequence of {\bf rational functions}, let's call them
$\{f^{(s)}_n(t)\}_{n=1}^{\infty}$, defined by
$$
f^{(s)}_n(t):=\sum_{k=0}^{n} \, \frac{n!^s}{(1+t)_k^s (1-t)_{n-k}^s} 
\quad,
$$
(note that $f^{(s)}_n(0)=F_n^{(s)}$). 

Their key idea was to apply the Zeilberger algorithm to the modified sum rather than the original sum, and see what happens. Let's recall some basic definition from
{\bf Wilf-Zeilberger algorithmic proof theory} [PWZ] that would also be needed later on, when we will describe our generalization of the Straub-Zudilin [SZ] work.

{\bf Definition 3} ([PWZ], p. 64): A discrete function $F(n,k)$ (defined on  $\{(n,k) \, | \, 0 \leq n,k <\infty \}$) is said to be a {\it proper hypergeometric term} if it
can be written in the form
$$
F(n,k)=P(n,k) \, 
\frac
{\prod_{i=1}^{uu} (a_i n + b_i k+ c_i)!}
{\prod_{i=1}^{vv} (u_i n + v_i k+ w_i)!}
x^k \quad,
\eqno(1)
$$

in which $x$ is an indeterminate over, say, the complex numbers, and

$\bullet$ $P$ is a polynomial

$\bullet$ the $a$'s, $b$'s, $u$'s, and $v$'s are specific integers, that is to say, they do not contain additional parameters, and

$\bullet$ the quantities $uu$ and $vv$ are finite non-negative, specific integers.

Recall that for any proper hypergeometric term $F(n,k)$ as defined above, the Zeilberger algorithm [Ze1] [PWZ] furnishes, for some non-negative integer $L$ (called the order) , polynomials
in $n$, $p_0(n), \dots , p_L(n)$, as well as another proper hypergeometric term, $G(n,k)$ called the {\bf certificate} (that furthermore has the property that
$G(n,k)/F(n,k)$ is a {\bf rational function} of $(n,k)$), such that
$$
p_0(n) F(n,k) \, + \,
p_1(n) F(n+1,k) \, + \, \dots \, +
p_L(n) F(n+L,k) \, = \,
G(n,k+1)-G(n,k) \quad .
$$

By summing with respect to $k$ from $k=0$ to $k=n$, we have the immediate corollary that the {\bf hypergeometric sum} 
$$
f(n):=\sum_{k=0}^{n} F(n,k) \quad ,
$$
satisfies the {\bf linear recurrence equation with polynomial coefficients}
$$
p_0(n) f(n) \, + \,
p_1(n) f(n+1) \, + \, \dots \, +
p_L(n) f(n+L) \, = \,
G(n,n+1)-G(n,0) \quad .
$$

Note that in general, the right side is not zero, so one gets an {\bf inhomogeneous} linear recurrence, but whenever the summand is {\bf natural}, for example
any {\bf binomial coefficient sum} that contains ${{n} \choose {k}}$ in it (in particular, of course, for the Franel sequences) the right side
vanishes, and one gets that the sum $f(n)$ satisfies a  {\bf homogeneous} linear recurrence equation with polynomial coefficients.

Since $(1+t)_k=(t+k)!/t!$ and $(1-t)_{n-k}=(n-(t+k))!/(-t)!$, the generalized Franel sum above, namely $f^{(s)}_n(t):=\sum_{k=0}^{n} \, \frac{n!^s}{(1+t)_k^s (1-t)_{n-k}^s} $,
yields the {\bf same recurrence}, but of course {\bf not} the same certificate. In particular the right side $G(n,n+1)-G(n,0)$  is no longer $0$ but some
multiple of the summand by a rational function of $t$ and $n$. But then came a {\bf nice surprise}, indeed, yet-another-{\it  miracle}.
For any {\it specific} $s$,  this miracle follows immediately from the output of the  Zeilberger algorithm. It is  proved {\it in general} in [SZ] (using {\bf human ingenuity}).

{\bf The Straub-Zudilin-Franel miracle}

For any $s>2$, the right side of the inhomogeneous recurrence satisfied by $f^{(s)}_n(t)$ is divisible by $t^{s+1}$ if $s$ is odd and by $t^s$ if $s$ is even.

It follows that for any $r<s$, the sequence obtained by extracting the coefficient of $t^r$ from the summand  $\frac{n!^s}{(1+t)_k^s (1-t)_{n-k}^s}$ of 
$f^{(s)}_n(t)$ automatically satisfies the {\bf very same} {\it homogeneous} linear recurrence satisfied by the Franel sequence $F^{(s)}_n$, but, of course, with different
initial conditions. Hence the limit of the ratios of the later sequence with the Franel sequence $F^{(s)}_n$ is an Ap\`ery limit of {\it some constant}.
Can we describe this constant directly? Yes we can!

Let's rewrite the sum $f^{(s)}_n(t)$ as follows

$$
f^{(s)}_n(t):=\sum_{k=0}^{n} \, \frac{n!^s}{(1+t)_k^s (1-t)_{n-k}^s } \, = \,
\sum_{k=0}^{n} \, \frac{n!^s}{k!^{s}(n-k)!^s} \cdot \left ( \frac{k!^s (n-k)!^s}{(1+t)_k^s (1-t)_{n-k}^s} \right ) \quad .
$$

Let's give the coefficient of $t^r$ in the Taylor expansion  (about $t=0$) of $\frac{k!^s (n-k)!^s}{(1+t)_k^s (1-t)_{n-k}^s}$ a name.

{\bf Definition 4}: $c^{(s)}_r(n,k)$ is the coefficient of $t^r$ in the Taylor expansion of the rational function of $t$  $\frac{k!^s (n-k)!^s}{(1+t)_k^s (1-t)_{n-k}^s}$.

{\bf Interesting consequences of the Straub-Zudilin-Franel miracle}

Now fix a power $s>2$ and $r<s+1$, and let $B(n)$ be the Franel sequence and A(n) the sequence of the coefficients of $t^r$ in the Taylor expansion of  $f^{(s)}_n(t)$. Then
$$
B(n)=\sum_{k=0}^{n} \, {{n} \choose {k}}^s \quad ,
$$
$$
A(n)=\sum_{k=0}^{n} \, {{n} \choose {k}}^s \, c^{(s)}_r(n,k) \quad .
$$
Now consider the Ap\'ery limit
$$
\lim_{n \rightarrow \infty} \frac{A(n)}{B(n)} =
\frac
{\sum_{k=0}^{n} \, {{n} \choose {k}}^s  c^{(s)}_r(n,k)}
{\sum_{k=0}^{n} \, {{n} \choose {k}}^s } \quad .
$$
Note that the left side is a certain {\it Ap\'ery limit} with the Franel recurrence and appropriate initial conditions (hence the limit can be computed very fast, with exponentially-decaying error). On the
other hand the right side is a {\it weighted-average}  of the $n+1$ numbers
$$
\{ c^{(s)}_r(n,k) \, | \, 0 \leq k \leq n \} \quad,
$$
with {\it weights}  $\{{{n} \choose {k}}^s\, |\, k=0 \dots n\}$. Since (recall the {\it Central Limit Theorem}) most of the weight is in the middle, we have that
the above Ap\'ery limit has an explicit description, namely
$$
\lim_{n \rightarrow \infty} c^{(s)}_r(n,\lfloor n/2 \rfloor) \quad.
$$

For any {\it specific} $r$, one can express $c^{(s)}_r(n,k)$  in terms of partial sums of harmonic-type series or series of type $\zeta(i)$. As $r$ gets larger, things get complicated
but our Maple package {\tt Zudilin.txt} can handle it easily (it is implemented by procedure {\tt CrL(A,r,L)}).
In particular, for $r=2$, things are still fairly easy. We have the following lemma, whose proof is omitted.

{\bf Lemma 1}:
$$
c^{(s)}_2(n,k)=\frac{s}{2} \left ( \, \sum_{i=1}^{k} \frac{1}{i^2} \, + \, \sum_{i=1}^{n-k} \frac{1}{i^2} \, \right ) \, + \,
\frac{s^2}{2} \left ( \, -\sum_{i=1}^{k} \frac{1}{i} \, +  \,  \sum_{i=1}^{n-k} \frac{1}{i} \right )^2 \quad.
$$

Note that in particular, the Ap\'ery limit for the case $r=2$ in the Straub-Zudilin $t$-version of Franel equals $s\,\zeta(2)$.

For each specific integer $s$, (and also for each specific even integer $r$) we can prove that (for sufficiently large $s$) you get $\zeta(r)$, as is proved, in general, using
human ingenuity, in [SZ].

{\bf Generalized Franel Sums}

The Zeilberger algorithm, just as easily, can find a homogeneous linear recurrence equation with polynomial coefficients for the generalized  Franel sum
$$
\sum_{k=0}^{n} \, {{n} \choose {k}}^s a^k \quad,
$$
for {\it any} positive rational number $a$ (or for that matter even for symbolic $a$). It turns out that for all the $s$ that we tried, the same miracle, i.e. that
the right side of the linear {\bf inhomogeneous} recurrence satisified by
$$
\sum_{k=0}^{n} \, \frac{n!^s}{(1+t)_k^s (1-t)_{n-k}^s} \cdot a^k \quad,
$$
is divisible by $t^r$ for $r \leq s$.  It follows that the first $s$ Taylor coefficients vanish. We are sure that this is provable in general using the method of [SZ],
but we leave it to  the interested reader.

Hence we have the following

{\bf Theorem 1}: For any positive rational number $a$, and for $s\geq 3$, and for $1 \leq r<s$ if $s$ is odd, and for $1 \leq r \leq s-1$ if $s$ is even.
Let $\alpha$ be the such that, as $n \rightarrow \infty$, $k=\alpha n$ maximizes the summand
$$
{{n} \choose {k}}^s a^k \quad
$$
(this can be easily found by computing the ratio of consecutive terms, setting it equal to $1$ and solving for $\alpha$ as $n \rightarrow \infty$, it is implemented in our Maple package
by procedure {\tt FindMaxk(F,n,k)}).
(Note that $\alpha$ is always an algebraic number, that our Maple package can always find in each case.)
Let $B(n)$ be $\sum_{k=0}^{n} \, {{n} \choose {k}}^s a^k$ and let $A(n)$ be solution of the same recurrence, of order $L$, say, that is satisfied by $B(n)$ by the Zeilberger algorithm, but
with the initial conditions extracted from the coefficients of $t^r$ in   
$
\sum_{k=0}^{n} \, \left ( \frac{n!}{(1+t)_{k} (1-t)_{n-k}} \right )^s a^k \quad,
$ 
for $n=0, \dots , n=L-1$, then the sequence $\frac{A(n)}{B(n)}$ converges {\bf with exponentially decaying error} to the constant

$$
\lim_{n \rightarrow \infty}   c^{(s)}_r(n,\lfloor \alpha n \rfloor) \quad .
$$

Note that trying to evaluate this limit from the definition would take for ever, since the convergence is so slow. So we can get lots of {\it minor Ap\'ery miracles}.

The output files {\tt https://sites.math.rutgers.edu/\~{}zeilberg/tokhniot/oZudilin1.txt} and
{\tt https://sites.math.rutgers.edu/\~{}zeilberg/tokhniot/oZudilin2.txt}  contain many examples.

{\bf The Zudilin-Straub t-Transform}

The beautiful construction in [SZ] was obtained by replacing $k!$ by $(1+t)_k$ and $(n-k)!$ by $(1-t)_{n-k}$. This naturally leads to

{\bf Definition 5}: The Zudilin-Straub $t$-transform of the proper hypegeometric term given in (1) is
$$
\hat{F}(n,k;t)=P(n,k+t) \, 
\frac
{\prod_{i=1}^{uu} ( b_i t+ 1)_{a_in+b_ik+c}}
{\prod_{i=1}^{vv}  ( v_i t+ 1)_ {u_i n + v_i k+ w_i}}
x^k \quad .
$$

It is immediate to see that the recurrence obtained via the Zeilberger algorithm applied to the Zudilin-Straub $t$-transform of any proper hypergeometric term is 
the same as the original, {\bf except} that the right side is not $zero$, i.e. the linear recurrence is {\it inhomogeneous}. Unfortunately, in general, the Franel-Straub-Zudilin
miracle does not occur, but it is easy to see that for any {\bf specific} positive integer $r$ the coefficient of $t^r$ in the Taylor expansion of the right side is
still $P$-recursive in $n$ (i.e. satisfies a linear recurrence equation with polynomial coefficients), and a recurrence for it can be algorithmically obtained, either by the `holonomic machine' ([ApaZ] [K]), or by `guessing' that can
be made fully rigorous using the general theorems of [ApaZ]. In this way we can get lots of {\it generalized Ap\'ery limits} describing constants defined
as limits of explicit sequences that converge very slowly. Rather than stating the formal theorem we refer the reader to the output files
given in the front of this article

{\tt https://sites.math.rutgers.edu/\~{}zeilberg/mamarim/mamarimhtml/zudilin.html} \quad .

{\bf The Maple package} {\tt Zudilin.txt} implements everything. Once you download it to your computer (that has Maple), load it to a Maple session with {\tt read `Zudilin.txt`}.
To get a list of the main functions, type {\tt ezra();}, and to get help with a specific function, type
{\tt ezra(FunctionName)}. For example, to get help with procedure {\tt ZT} (that implements the Zudilin-Starub $t$-transform) type {\tt ezra(ZT);}.
Enjoy!

{\bf Conclusion}

While major Ap\'ery miracles are few and far between, thanks
to the {\bf Zudilin-Straub} $t$-transform, we can obtain many minor Ap\'ery miracles. The front of this article contains numerous examples, but using
our Maple package, readers can generate many new ones.

\bigskip

{\bf References}

[A] Roger Ap\'ery, 
{\it ``Interpolation de fractions continues et irrationalit\'e
de certaine constantes''}
Bulletin de la section des sciences du C.T.H.S. \#3
p. 37-53, 1981.

[ApaZ] Moa Apagodu  and Doron Zeilberger,
{\it Multi-variable Zeilberger and Almkvist-Zeilberger algorithms and the
sharpening of Wilf-Zeilberger Theory },
Adv. Appl. Math. {\bf 37} (2006)(Special Regev issue), 139-152. \hfill\break
{\tt https://sites.math.rutgers.edu/\~{}zeilberg/mamarim/mamarimhtml/multiZ.html} \quad .

[CS] Marc Chamberland and Armin Straub, {\it Ap\'ery limits: experiments and proofs},   
 Amer. Math. Monthly {\bf 128} (2021), 811–824.  \hfill\break
{\tt https://arxiv.org/abs/2011.03400} \quad .

[DKZ] Robert Dougherty-Bliss, Christoph Koutschan, and Doron Zeilberger,
{\it Tweaking the Beukers integrals in search of more miraculous irrationality proofs \'a la Ap\'ery},
Ramanujan Journal, published-on-line before print, volume and page number tbd. \hfill\break
{\tt https://sites.math.rutgers.edu/\~{}zeilberg/mamarim/mamarimhtml/beukers.html}

[DZ] Robert Dougherty-Bliss and Doron Zeilberger, {\it Experimenting with Ap\'ery limits and WZ pairs},
Maple Transactions, {\bf v.1}, issue 2 (2021). \hfill\break
{\tt https://mapletransactions.org/index.php/maple/article/view/14359} \hfill\break
{\tt https://sites.math.rutgers.edu/\~{}zeilberg/mamarim/mamarimhtml/wzp.html} \quad .

[K] Christoph Koutschan, {\it   Advanced applications of the holonomic systems approach}, 
PhD thesis, Research Institute for Symbolic Computation (RISC), Johannes Kepler University, Linz, Austria, 2009.\hfill\break
{\tt http://www.koutschan.de/publ/Koutschan09/thesisKoutschan.pdf} \quad .

[PWZ] M. Petkovsek, H. S. Wilf and D. Zeilberger,
{\it A=B}, AK Peters, Wellesley, (1996). \hfill\break
{\tt https://sites.math.rutgers.edu/\~{}zeilberg/AeqB.pdf} \quad .

[SZ] Armin Straub and Wadim Zudilin, {\it Sums of powers of binomials, their Ap\'ery limits, and Franel suspicions},
Intern. Math. Research Notices (to appear), 19 pages. \hfill\break
{\tt https://arxiv.org/abs/2112.09576} \quad .

[Ze1] Doron Zeilberger, {\it The method of creative telescoping}, 
J. Symbolic Computation {\bf 11} (1991), 195-204.  \hfill\break
{\tt https://sites.math.rutgers.edu/\~{}zeilberg/mamarimY/creativeT.pdf} \quad .

[Ze2] Doron Zeilberger, {\it Closed form (pun intended!)}, 
 in: ``Special volume in memory of Emil Grosswald", M. Knopp
and M. Sheingorn,
eds., Contemporary Mathematics {\bf 143}, 579-607, AMS, Providence (1993). \hfill\break
{\tt https://sites.math.rutgers.edu/\~{}zeilberg/mamarim/mamarimhtml/pun.html} \quad ,

[Ze3] Doron Zeilberger, {\it Computerized deconstruction}, Advances in Applied Mathematics {\bf 30} (2003), 633-654. \hfill\break
{\tt https://sites.math.rutgers.edu/\~{}zeilberg/mamarim/mamarimhtml/derrida.html} \quad .

[ZeZu] Doron Zeilberger and Wadim Zudilin, {\it Automatic Discovery of Irrationality Proofs and Irrationality Measures},
International Journal of Number Theory {\bf 17} (2021), 815-825. \hfill\break
{\tt https://sites.math.rutgers.edu/\~{}zeilberg/mamarim/mamarimhtml/gat.html}

\bigskip
\hrule
\bigskip
Robert Dougherty-Bliss , Department of Mathematics, Rutgers University (New Brunswick), Hill Center-Busch Campus, 110 Frelinghuysen
Rd., Piscataway, NJ 08854-8019, USA. \hfill\break
Email: {\tt  robert dot w dot bliss at gmail dot com}   \quad .
\bigskip
Doron Zeilberger, Department of Mathematics, Rutgers University (New Brunswick), Hill Center-Busch Campus, 110 Frelinghuysen
Rd., Piscataway, NJ 08854-8019, USA. \hfill\break
Email: {\tt DoronZeil at gmail  dot com}   \quad .

\end